\documentclass[12pt]{article}
 
\usepackage{a4,amsmath,amssymb,amsthm,url}
\usepackage{verbatim}

\usepackage{tikz}
\usetikzlibrary{decorations.pathreplacing} 
\usetikzlibrary{shapes}
\usetikzlibrary{arrows}

\newtheorem{theorem}{Theorem}

\newcommand{\N}{{\mathbb N}}
\newcommand{\F}{{\mathcal F}}
\newcommand{\A}{{\mathcal A}}
\newcommand{\K}{{\mathcal K}}

\newcommand{\R}{{\mathbb R}}

\newcommand\dsum{\sum\limits}

\theoremstyle{plain}

\theoremstyle{definition}

\theoremstyle{remark}
\newtheorem{example}{Example}

\newcommand {\Bild} {\operatorname{im}}
\newcommand {\Hom}{\operatorname{Hom}}

\newcommand{\abs}[1]{\lvert#1\rvert}

\newcommand\dateiname[1]{}

\newcommand\datei[1]{#1}

\renewcommand\bar[1]{\overline{#1}}

\title{Counting monochromatic copies of $K_4$: a new lower bound for the Ramsey multiplicity problem}
\author{Susanne Nie{\ss}\footnote{
The author was supported by DFG project TA 309/2-2 and a scholarship of the state of Bavaria according to the Bayerisches Elitef\"orderungsgesetz (BayEFG) 
}
\\
	TU M\"{u}nchen}
\date{}

\begin{document} 
\maketitle

\begin{abstract}
Denote by $k_4(n)$ the minimal number of monochromatic copies of a $K_4$ in a 2-colouring of the edges of $K_n$ 
and let $c_4 :=\lim_{n\rightarrow\infty} {k_4(n)} / {\binom {n}{4}}$.  
The best known bounds so far were given by Thomason, who proved that
$c_4<\frac { 1 } { 33 }\approx 0.0303$, 
and Giraud, who showed that 
$c_4 > \frac { 1 } { 46 }\approx 0.0217$.  
In this paper we prove the new lower bound $c_4 > \frac{204603019}{7112448000} > 0.0287$.
\end{abstract}

\section{Introduction}
\label{sec:intro}
We denote by $ k_t(G) $ the number of cliques on $t$ vertices in a graph $G$ and define 
$k_t(n):=\min\{ k_t(G)+k_t(\bar { G }):\abs{ G } = n\}$. In other words, $k_t(n)$ is the minimum number of 
monochromatic copies of $K_t$ in a $2$-edge-colouring of $K_n$. 
It follows from Ramsey's theorem that $k_t(n)>0$, if $n$ is sufficiently large compared to $t$. We now let
$$ 
c_t :=\lim_{ n\rightarrow\infty }\frac { k_t(n)} {\binom { n } { t } }. 
$$

The problem of determining $ c_t $ is known as the Ramsey Multiplicity Problem
and was initiated by Erd\H{o}s~\cite{Erd}.  
For $ t = 3 $ it was shown by Goodman~\cite{Goo} that 
$ c_3 =\frac { 1 } { 4 } $,
but for $ t\geq 4 $ the problem of determining $c_t$ is still open. 
Currently, the best general bounds are as follows: 
Thomason proved in ~\cite{Tho} that $c_t \le 0,936 \cdot 2^{ 1 -\binom { t } { 2 } } $
On the other hand, not much work seems to have been done for the lower bounds, 
but the current record was only recently set by Conlon~\cite{Con}, who showed that
$ c_t\geq C ^ {-t ^ 2 \cdot(1+o(1)} $, for a constant $ C\approx 2.18 $.

For the case $t=4$, 
the best known upper bound for $c_4$ was given by Thomason~\cite{Tho} and states that 
$ c_4 <\frac { 1 } { 33 }\approx 0.0303 $, while the best known lower bound 
$ c_4 >\frac { 1 } { 46 }\approx 0.0217 $ so far
was due to Giraud~\cite{Gir}.\footnote{
Giraud also showed that if one could prove that
$k_4(n)$ is attained by a graph where the portion of $ K_3 $ and $\bar { K }_3 $ 
is smaller than $\frac { 1 } { 4 } $, then this would imply that 
$\frac { k_4(n) } {\binom { n } { 4 } }\geq\frac { 1 } { 35 }\approx 0,02857 $.}
 In this short paper we prove the following new lower bound:

\begin{theorem}
$$
c_4\geq \frac{204603019}{7112448000} >0.02876689
$$
\end{theorem}

The proof if this lower bound relies heavily on the concept of flag algebras introduced by 
Razborov~\cite{Raz}. This paper is organized as follows:  At first we give a brief introduction into flag agebras. Then we describe how we used them to prove theorem 1. After that we explain the data files that are given for checking the proof.

\section{The proof}

We first introduce some notation. For a graph $ G $, denote by 
$ V(G)$ the vertex set of $ G $ and by $ E(G)$ its edge set and 
let $ v(G)=\abs { V(G) } $ und $ e(G)=\abs {E(G)} $. 
For a natural number $n$ we set $[n]= \{ 1,\ldots, n\} $. 

A graph $\sigma $ with vertex set $ V(\sigma)=[s]$ is called a type. 
A pair $ F =(G,\theta)$ is called a $\sigma $-flag, if $ G $ is a graph and 
$\theta:[s]\to V(G)$ an injective function such that 
$ \{\theta(i),\theta(j)\}\in E(G)\Leftrightarrow \{ i, j\}\in E(\sigma) $. 
In particular, every graph can be considered as a 0-flag with $ s = 0 $ and 
$\sigma =\varnothing $. In the following, we always denote by $s$ the number of vertices in 
$\sigma $. 

For a flag $ F =(G,\theta)$ we let $ V(F):= V(G)$. 
We call two $\sigma $-flags $ F_1 =(G_1,\theta_1)$ and $ F_2 =(G_2,\theta_2)$ isomorphic and 
write $ F_1\cong F_2 $, 
if there is a bijective mapping $\psi: V(G_1)\to V(G_2) $ such that 
$\psi\theta_1 =\theta_2 $ and 
$$
\forall v, w\in V(G_1):(\{\psi(v),\psi(w)\}\in E(G_2)\Leftrightarrow\{ v, w\}\in E(G_1)) .
$$
For a type $\sigma $ of order $ s\in\N_0 $ und $ \ell\geq s $ let 
$\F ^ {\sigma }_\ell $ be the set of all $\ell$-vertex $\sigma $-flags up to isomorphism 
and set $\F ^ {\sigma } :=\bigcup _ { \ell = s } ^\infty\F ^ {\sigma } _\ell $. 
For a flag $ F = (G,\theta) $ and a set $ U\subseteq V (F) $ with 
$\Bild(\theta)\subseteq U $ let $F\mid_U $ be the flag obtained from $ F $
if one replaces $ G $ by $G[U]$, the subgraph of $G$ induced by the set $U$.

In a flag $ F = (G,\theta) $ a finite family of subsets $ V _ i\subseteq V (G) $ 
with the property $ V _ i\cap V _ j =\Bild(\theta) $ for $i\neq j$ is called a sunflower and 
the sets $ V_i $ are called the petals of the sunflower.
For $ n\geq 1 $ and $\sigma $-flags $ F _ 1,\ldots, F _ n $ and a $\sigma $-flag 
$ F $ with $\dsum_{ i = 0 }^n (v (F_i) - s)\leq v (F) - s $, we define the 
density $ p (F_1,\ldots, F_n; F) $ as follows:
Let $B$ be the set of all sunflowers $ (V_1,\ldots, V_n) $ in $ V (F) $ with 
$\abs { V_i } = v (F_i)$ for all $i\in [n] $, then 
$$
p (F_1,\ldots, F_n; F) =\frac {\abs {\{ (V_1,\ldots, V_n)\in B\colon \quad
F\mid_{V_i}\cong F_i\quad\forall i\in [n]\} } } {\abs { B } } .
$$

\begin{figure}[ht]
\begin{center}
    \begin{tikzpicture} 

\draw (1,5) node (a1) [] { };
\draw (2,5) node (b1) [] { };
\node at (1.5, 4.5) { $\rho$ };
    \draw (a1) -- (b1);    

\draw (3,5) node (a2) [] { };
\draw (4,5) node (b2) [] { };
\node at (3.5,4.5) { $\bar {\rho } $ };

\draw (5,5) node (a3) [] { };
\draw (6,5) node (b3) [] { };
\draw (5.5,6) node (c3) [] { };
\draw (a3) -- (b3) -- (c3) -- (a3);
\node at (5.5, 4.5) { $ K_3 $ };

\draw (7,5) node (a4) [] { };
\draw (8,5) node (b4) [] { };
\draw (7,6) node (c4) [] { };
\draw (8,6) node (d4) [] { };
\draw (a4) -- (b4) -- (d4) -- (c4) -- (a4);
\node at (7.5,4.5) { $ G_1 $ };

\draw (9,5) node (a5) [] { };
\draw (10,5) node (b5) [] { };
\draw (9,6) node (c5) [] { };
\draw (10,6) node (d5) [] { };
\draw (a5) -- (b5) -- (d5) -- (c5) -- (a5);\draw (a5) -- (d5);
\node at (9.5,4.5) { $ G_2 $ };

\draw (1,3) node (a6) [label=270:{\tiny 1}] { };
\draw (2,3) node (b6) [] { };
\draw (a6) -- (b6);
\node at (1.5,2.5) { $ e $ };

\draw (3,3) node (a7) [label=270:{\tiny 1}] { };
\draw (4,3) node (b7) [] { };
\node at (3.5,2.5) { $\bar { e } $ };

\draw (5,3) node (a8) [label=270:{\tiny 1}] { };
\draw (6,3) node (b8) [] { };
\draw (5.5,4) node (c8) [] { };
\draw (a8) -- (b8) -- (c8) -- (a8);
\node at (5.5,2.5) { $ K_3^1 $ };

\draw (7,3) node (a9) [label=270:{\tiny 1}] { };
\draw (8,3) node (b9) [] { };
\draw (7.5,4) node (c9) [] { };
\draw (a9) -- (b9) -- (c9);
\node at (7.5, 2.5) { $ P_3^{ 1, b } $ };

\draw (9,3) node (a10) [label=270:{\tiny 1}] { };
\draw (10,3) node (b10) [] { };
\draw (9.5,4) node (c10) [] { };
\draw (b10) -- (a10) -- (c10);
\node at (9.5,2.5) { $ P_3^{ 1, c } $ };

\draw (11,3) node (a11) [label=270:{\tiny 1}] { };
\draw (12,3) node (b11) [] { };
\draw (12,4) node (c11) [] { };
\draw (11,4) node (d11) [] { };
\draw (a11) -- (b11) -- (c11) -- (d11) -- (a11);
\node at (11.5,2.5) { $ H_1 $ };

\draw (13,3) node (a12) [label=270:{\tiny 1}] { };
\draw (14,3) node (b12) [] { };
\draw (14,4) node (c12) [] { };
\draw (13,4) node (d12) [] { };
\draw (a12) -- (b12) -- (c12) -- (d12) -- (a12);\draw (a12) -- (c12);
\node at (13.5,2.5) { $ H_2 $ };

\draw (15,3) node (a13) [label=270:{\tiny 1}] { };
\draw (16,3) node (b13) [] { };
\draw (16,4) node (c13) [] { };
\draw (15,4) node (d13) [] { };
\draw (a13) -- (b13) -- (c13) -- (d13) -- (a13);\draw (b13) -- (d13);
\node at (15.5,2.5) { $ H_3 $ };

\draw (1,1) node (a14) [label=270:{\tiny 1}] { };
\draw (2,1) node (b14) [label=270:{\tiny 2}] { };
\draw (a14) -- (b14);
\node at (1.5,0.5) { $\sigma_1 $ };

\draw (3,1) node (a15) [label=270:{\tiny 1}] { };
\draw (4,1) node (b15) [label=270:{\tiny 2}] { };
\draw (3.5,2) node (c15) [] { };
\draw (b15) -- (a15) -- (c15);
\node at (3.5,0.5) { $ H_4 $ };

\draw (5,1) node (a16) [label=270:{\tiny 1}] { };
\draw (6,1) node (b16) [label=270:{\tiny 2}] { };
\draw (5.5,2) node (c16) [] { };
\draw (a16) -- (b16) -- (c16);
\node at (5.5,0.5) { $ H_5 $ };

\draw (7,1) node (a17) [label=270:{\tiny 1}] { };
\draw (8,1) node (b17) [label=270:{\tiny 2}] { };
\draw (8,2) node (c17) [] { };
\draw (7,2) node (d17) [] { };
\draw (a17) -- (b17) -- (c17) -- (d17) -- (a17);
\node at (7.5,0.5) { $ H_6 $ };

\draw (11,1) node (a18) [label=270:{\tiny 1}] { };
\draw (12,1) node (b18) [label=270:{\tiny 2}] { };
\node at (11.5,0.5) { $\sigma_2 $ };

\draw (13,1) node (a19) [label=270:{\tiny 1}] { };
\draw (14,1) node (b19) [label=270:{\tiny 2}] { };
\draw (13.5,2) node (c19) [] { };
\draw (a19) -- (c19) -- (b19);
\node at (13.5,0.5) { $ H_7 $ };

\draw (15,1) node (a20) [label=270:{\tiny 1}] { };
\draw (16,1) node (b20) [label=270:{\tiny 2}] { };
\draw (16,2) node (c20) [] { };
\draw (15,2) node (d20) [] { };
\draw (a20) -- (c20) -- (b20) -- (d20) -- (a20);
\node at (15.5,0.5) { $ H_8 $ };

   \foreach \a in {a,b} {
    	\foreach \i in {1,2,6,7, 14, 18} {
    		\fill [black,opacity=.99] (\a\i) circle (3pt);
    	}
    }
   \foreach \a in {a,b,c} {
    	\foreach \i in {3,8,9,10, 15, 16, 19} {
    		\fill [black,opacity=.99] (\a\i) circle (3pt);
    	}
    }
   \foreach \a in {a,b,c,d} {
    	\foreach \i in {4,5, 11, 12,13, 17, 20} {
    		\fill [black,opacity=.99] (\a\i) circle (3pt);
    	}
    }
    \end{tikzpicture}
\caption{Some examples for flags}
\label{fig:graph1}
\end{center}
\end{figure}
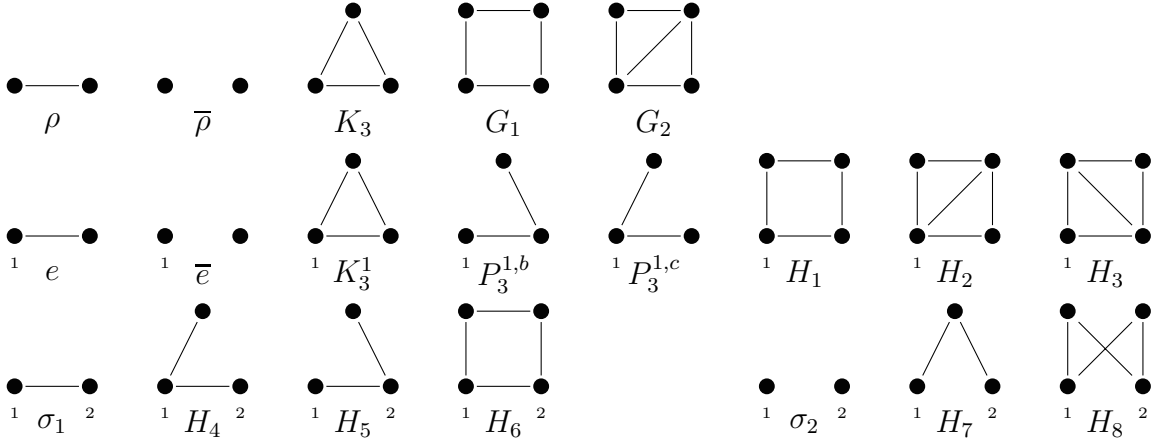

\begin{example}
In figure 1, we give a few examples. Here $\rho,\ \bar {\rho },\ K_3,\ G_1 $ and $ G_2 $ are 0-flags, that is graphs, with $ p(\rho; K_3)=\frac { 3 } { 3 } = 1,\ p(\rho; G_1)=\frac { 4 } { 6 }=\frac { 2 } { 3 },\ p(\rho; G_2)=\frac { 5 } { 6 },\ p(\bar {\rho }; K_3)= 0,\ p(\bar {\rho }; G_2)=\frac { 1 } { 6 },\ p(K_3; G_1)= 0,\ p(K_3; G_2)=\frac { 2 } { 4 } =\frac{1}{2} $. For $s=0$, a sunflower with 2 petals consists of 2 disjoint subsets, so we have $ p(\rho,\rho; G_1)=\frac { 4 } { 6 },\ p(\bar {\rho },\bar {\rho }; G_1)=\frac { 2 } { 6 },\ p(\rho,\bar {\rho }; G_1)= 0,\ p(\rho,\rho; G_2)=\frac { 4 } { 6 },\ p(\rho,\bar {\rho }; G_2)= p(\bar {\rho },\rho; G_2)=\frac { 1 } { 6 },\ p(\bar {\rho },\bar {\rho }; G_2)= 0 $. For $ s=1 $, the only type is an isolated vertex. In the example above $ e,\ \bar { e },\ K_3 ^ 1,\ P_3 ^ { 1, b },\ P_3 ^ { 1, c },\ H_1,\ H_2,\ H_3 $ one can see that sometimes there are different flags corresponding to one graph, for example $ H_2 $ and $ H_3 $ both correspond to $ G_2 $. Here we have $ p(e; K_3 ^ 1)= p(e; P_3 ^ { 1, c })= \frac { 2 } { 2 } = 1 $ but $ p(e, P_3 ^ { 1, b } =\frac{1}{2} $ and $ p(e, H_1)= p(e, H_3)=\frac { 2 } { 3 } $ but $ p(e, H_2)=\frac { 3 } { 3 } $. On the other hand, $ p(K_3 ^ 1; H_1)= 0,\ p(K_3 ^ 1 H_2)=\frac { 2 } { 3 },\ p(K_3 ^ 1; H_3)=\frac { 1 } { 3 } $ and finally $ p(\bar { e }, K_3 ^ 1; H_3)=\frac { 1 } { 3 },\ p(e,\bar { e }; H_1)= p(\bar { e }, e; H_1)=\frac { 2 } { 6 } $. For $s=2$, there are the two types $\sigma_1 $ and $\sigma_2 $. With the $\sigma_1$-flags $H_4,\ H_5,\ H_6$ we have $ p(H_4; H_6)= \frac{1}{2},\ p(H_4, H_4; H_6)= 0 $ and $ p(H_4, H_5; H_6)= \frac{1}{2} $. With the $\sigma_2$-flags $H_7$ and $H_8$ we have $ p(H_7; H_8)=\frac { 2 } { 2 } = 1 $.
\end{example}

Denote by $\R\F^{\sigma } $ the space of all finite formal sums of $\sigma $-flags 
with coefficients in $\R $. We write $\K^{\sigma }$ for the subspace generated by all elements of the form 
$$
\tilde { F } -\dsum_{ F \in\F^\sigma_{ \ell_2 } } p (\tilde { F }; F)F  
$$ 
with $ F\in\F^{\sigma }_{ \ell_1 }, $ $ \ell_1\leq \ell_2 $ 
and let $\A^\sigma :=\R\F^\sigma/\K^\sigma $. 
For two $\sigma $-flags $ F_1\in\F^\sigma_{ l_1 }, F_2\in\F^\sigma_{ \ell_2 } $ and 
$ \ell\geq \ell_1 + \ell_2 - s $ let 
$$ 
F_1\cdot F_2 :=\dsum_{ F\in\F^\sigma_\ell } p (F_1, F_2; F) F. 
$$
Owing to the structure of $\A^\sigma $, this product is independent of the choice of 
$ \ell $ and thus well defined. 
By linear extension of the product,  $\A^\sigma $ forms an algebra with the neutral element 
$\sigma $. 

For $\sigma $ on $ s >0 $ vertices and a $\sigma $-flag $F = (G,\theta)$ 
we define  $q_\sigma (F) $ as follows: 
Let $ \Psi$ be the set of all injective functions $\psi: [s]\rightarrow V (G) $ and set
$$
q_\sigma (F) 
=\frac {\abs {\{\psi\in \Psi\colon G\mid_{\Bild (\psi) } =\sigma\wedge (G,\psi)\cong F\} } } 
{\abs { \Psi } } .
$$
Then we can construct a linear function $ [[. ]]_\sigma:\A^{\sigma }\rightarrow\A^0 $ 
by mapping every $\sigma $-flag $ F = (G,\theta) $ to the value 
$ [[ F ]]_\sigma = q_\sigma (F)\cdot G $. 

Let $\Hom^+ (\A^\sigma) $ denote the set of all algebra homomorphisms 
$\phi:\A^\sigma\rightarrow\R $ with $\phi (F)\geq 0$ for every $F\in\F^\sigma $. 
Razborov~\cite{Raz} showed that every $\phi\in\Hom^+ (\A^\sigma) $ 
can be written as the limit of a sequence of functions 
$ p (.; F_n) $ with $ v (F_n)\rightarrow\infty $,
and that, vice versa, every limit of such a convergent sequence lies in $\Hom^+ (\A^\sigma) $. 

Moreover, for $ n\in\N $ and a vector $ x\in (\F^\sigma)^{ n\times 1 } $ 
whose components are $\sigma $-flags und a symmetric, positive semidefinite matrix 
$ A\in\R^{ n\times n } $, we know that 
$$
\phi (x^{T} A x)\geq 0\qquad \forall\phi\in\Hom^+ (\A^\sigma). 
$$
Here, for $ f, g\in\A^{\sigma } $ the inequality $ f\leq g $ means that 
$\phi (f)\leq\phi (g)$ for every $\phi\in\Hom^+ (\A^\sigma) $.

More details and proofs can be found in \cite{Raz}.

With the notation introduced, we can now reformulate the Ramsey multiplicity problem:  
$ k_4 (n) $ is the minimum number of $ K_4 $ and \emph{induced} $\bar { K }_4 $ 
in a graph with $ n $ vertices, $ c_4 (n) =\min\{ p (K_4, G) + p (\bar { K }_4, G)\mid v (G) = n\} $ 
and $ c_4 =\min\{\phi (K_4 +\bar { K }_4)\mid\phi\in\Hom^+ (\A^0)\} $. 

Hence we can proceed as follows: We choose types $\sigma_1,\ldots,\sigma_n $ und $ \ell\in\N $, 
construct vectors $ x_i $ from lists of $ F\in\F^{\sigma_i }_\ell $ and 
find positive semidefinite matrices $ M_1,\ldots, M_n $ and a real number $ c >0 $, 
so that for every $\phi\in\Hom^+ (\A ^0)$ 
\begin{equation}
\label{eq:grZero}
\phi (K_4 +\bar { K }_4) -\phi ([[ x_1^{ T } M_1 x_1 ]]_{\sigma_1 } +\ldots + 
[[ x_n^{ T } M_nx_n ]]_{\sigma_n }) - c \geq 0 
\end{equation}
Since all $\phi\in\Hom^+ (\A^0) $ are non-negative on all quadratic forms 
$ [[ x_i^{ T } M_i x_i ]]_{\sigma_i } $, it follows that $\phi (K_4 +\bar { K }_4)\geq c $. 
Since every flag in a flag-algebra can be written as the weighted sum of flags with a larger number of vertices
and since, by definition, every $\phi\in\Hom^+ (\A^0) $ is non-negative on every flag, 
we can check (\ref{eq:grZero}) by comparing coefficients. 

Choosing $s>0$ and $\sigma_1,\ldots,\sigma_n $ with $\abs {\sigma_i } = s $ as well as $ \ell_1 >s $ and constructing the vectors $x_i$ from lists of the $ F\in\F^{\sigma_i }_{ \ell_1 }$, one can write the quadratic forms $ x_i^T M_i x_i $ as weighted sums of flags in $\F^{\sigma_i }_{ \ell_2 } $ with $ \ell_2 := 2\cdot \ell_1 - s $. The multiplication of flags can be described by a vector of $ \abs {\F^{\sigma_i }_{ \ell_1 }} \times \abs {\F^{\sigma_i }_{ \ell_1 }} $-matrices $ D_j $ with $ 1\leq j\leq\abs {\F^{\sigma_i }_{ \ell_2 }} $, where the $ (j_1, j_2) $-component in $ D_j $ is $ p (F_{ j_1 }, F_{ j_2 }; F_j) $ with $ F_{ j_1 }, F_{ j_2 }\in\F^{\sigma_i }_{ \ell_1 }, $ $ F_j\in\F^{\sigma_i }_{ \ell_2 } $, the coefficient of $ F_j $ in the mulitplication of $ F_{ j_1 } $ with $ F_{ j_2 } $. Then  the $j$th coefficient of the quadratic form is achieved by multiplying $D_j$ with $M_i$ component-wise and then taking the sum of all components. The operator $ [[. ]]_{\sigma_i } $ turns this into weighted sums of flags in $\F^0_{ \ell_2 } $. Their coefficient vector $v_i$ is obtained by multiplying the current vector by the matrix of $ q_{\sigma } $. As for all relevant combinations of $s$ and $\ell_1$ always $\ell_2>4$, for $ K_4 +\bar { K }_4 $ the vector $w$ of the $ p (K_4; F_k) + p (\bar { K }_4; F_k) $ with $ F_k\in\F^0_{ \ell_2 } $ has to be used. Defining at last $e$ as the vector of lenght $\abs{\F^0_{ \ell_2 } } $ where each entry equals 1, the problem can be formulated as the search for positive definite matrices $M_i$ and the biggest possible $c>0$ so that $ v_1 +\ldots + v_n + c\cdot e\leq w $ component-wise.

Suitable matrices $ M_i $ can be found through semidefinite programming (SDP). 
For this, we have used {SeDuMi}~\cite{SeD} which runs on $Matlab$ or $Octave$. 
However, this yields matrices composed of floating point numbers, which cannot be used for a proof. Therefore,
we transformed these matrices with the help of $Maple$ into matrices containing only rational numbers, which 
are again positive semidefinite. To check this property, we computed the smallest eigenvalue (these are now the 
only floating point numbers left in the proof, all other numbers are integers or rationals). 

We have worked with $ s=3 $ and 4 different types: 
$\sigma_1 $ is the empty graph on vertex set $ [3] $ and the adjacency lists of the other types are 
$$\sigma_2: (1,2) $$
$$\sigma_3: (1,2), (1,3) $$
$$\sigma_4: (1,2), (1,3), (2,3). $$
With each of these $ \sigma $ we have constructed the quadratic forms on $\sigma $-flags with 5 vertices 
and represented these products as flags with 7 vertices. 
The results of the 
$ [[ x^{T} M x ]]_\sigma $ are then weighted sums of graphs with 7 vertices.
On the other hand, we have represented $ K_4 $ and $\bar { K }_4 $ also as sums of 7-vertex graphs and then 
compared their coefficients. 

In principle, this proof could be checked by hand and written out completely on paper, 
but due to the size of the matrices this would be very tiresome. Instead, we follow the approach of 
Pikhurko in~\cite{Pik} who suggested that all the data needed to check the proof should be made available in electronic form. It can be downloaded from

\url{http://www-m9.ma.tum.de/Allgemeines/SusanneNiess}

Here is a rough description of the content of the files that can be found there. 
Since the rational numbers $ p (F_1, F) $ are fractions which, for fixed number of vertices in 
$ F_1 $ and $ F $, all have the same denominator, we have only stored the numerator for the values of $p$,
and similarly for the values $ q_\sigma $.

Since the matrices can only be understood if the order of the flags is known, we first list the files containing the flags that we used. The files called \datei{jba35\_}$i$ with $ i\in [4] $ 
contain the 5-vertex $\sigma_i $-flags in $nauty$-format (see~\cite{McK}), while
the files \datei{jbc35\_}$i$ contain the same flags in a format that is more readable by humans: 
in the first line, there is an integer which denotes the number of flags, then, after this, 
each line contains the upper-triangular adjancey matrix of a flag, all rows one after the other. 

The 7-vertex $\sigma_i $-flags in $nauty$-format are stored in the files called \datei{jba37\_}$i$ with $ i\in [4] $, 
again mirrored in the files \datei{jbc37\_}$i$ as upper-triangular matrices. 
The file \datei{jba07} contains all 7-vertex 0-flags in $nauty$-format, the file \datei{jbc07} is again the more readble version.

The vector with the $ p (K_4, F) + p (\bar { K }_4, F) $ for all $ F\in\F^0_7 $, 
in other words the coefficients of $ K_4 +\bar { K }_4 $ represented as weighted sums of 7-vertex flags
is stored in the file \datei{l47}.
Here again we only saved the numerator, the common denominator is $\binom { 7 } { 4 } = 35 $.

Information about the factors $ q_{\sigma_i } (F) $ can be found in the files \datei{qjb37\_}$i$.
The first column contains the indices of the 
$ F\in\F^{\sigma_i } $, the second column the indices of the corresponding 0-flags, 
the third column the numerator of the $ q_{\sigma_i } (F) $. 
The denominator is 210 for all of them. 

The coefficients for the multiplication of the 5-vertex $\sigma_i$-flags can be found in the 
files \datei{si35\_}$i$. 
Here a column of the form ``$a\ b\ c\ d$'' means that 
the coefficient of the $a$-th 7-vertex $\sigma_i$-flag 
when multiplying the $b$-th and the $c$-th $5$-vertex $\sigma_i$-flag has the numerator $d$. 
The denominator is always $\binom{7-3}{5-3}=6$. 

The files \datei{si35\_}$i$ and \datei{qjb37\_}$i$ are used to compute the files \datei{sp35\_}$i$ which are 
in the same format as the \datei{si35\_}$i$ and contain the coefficients 
of the 7-vertex 0-flags in $[[F_b\cdot F_c]]_{\sigma_i}$. Here the denominator is always $210\cdot 6=1260$. 

The matrices $M_i$ describing the quadratic forms are stored in the files \datei{yy}$i$.\datei{m}, 
using the internal format of $Maple$. 
In addition, the numerators of the rationals in these matrices are saved in the files 
\datei{yz}$i$ in csv-format, the corresponding denominator is 11289600 and stored in the file \datei{yzn}. 

For these matrices it is true that: 
$$ 
[[x_1^{-1}M_1x_1]]_{\sigma_1} + [[x_2^{-1}M_2x_2]]_{\sigma_2} + 
[[x_3^{-1}M_3x_3]]_{\sigma_3} + [[x_4^{-1}M_4x_4]]_{\sigma_4} +\frac{204603019}{7112448000}\leq K_4+\bar{K}_4 
$$
where each $x_i$ is the vector of the 5-vertex $\sigma_i$-flags. 
By the nonnegativity of the quadratic form it follows that 
$ \frac{204603019}{7112448000}\leq K_4+\bar{K}_4 $ and thus $c_4\geq \frac{204603019}{7112448000}$.

For this proof, numerous computational difficulties had to be overcome. 
In order to push the lower bound for $c_4$ as far as possible, 
we had to produce lists of flags and matrices of the corresponding coefficient for $\ell$ as large as possible. As both the number of graphs and the number of $\sigma$-flags for any $\sigma$ grows hyperexponentially in $\ell$, this leads of course to problems with time, memory and disk space. In addition to improving computer hardware, it was necessary to optimize the efficiency of the 
software used. 
Both for listing the flags and for calculating many of the coefficients, $nauty$ by Brendan McKay~\cite{McK} led to great improvements. 
For calculating the coefficients $ p (F_1; F) $ and $ p (F_1, F_2; F) $, the algorithm 
in~\cite{RW} was useful for listing all sunflowers more efficiently. For some calculations we have written computer programmes where several functions are available for the same task, but with different time-memory-tradeoffs. 
For calculating the $q_{\sigma}$ it turned out  that the best method was to store some data in a raw form while producing the $\sigma$-flags and then processing them with $Matlab$ or $Octave$ ~\cite{Oct}. 
Using those programmes it turned out that $Octave$, an open source version of $Matlab$, was more suitable because it can process bigger amounts of data,
and therefore most of the necessary matrix-operations were done with $Octave$. 
While working on this problem, we learned about the development of $Flagmatic$ ~\cite{Vau}.
However, due to its file format $Flagmatic$ is confined to graphs with at most 9 vertices, and so we decided to stick to our own software which can cope with flags on 10 and 11 vertices (and it should 
be possible to increase this still further). 
Unfortunately, $SeDuMi$ is not able to process the matrices arising from flags with such a high number of vertices because of their size but they are still available for experiments with other SDP-solvers. Among the matrices that we can process up to now, the ones described above produced the highest lower bound for $c_4$.

\section{Acknowledgements}
The author would like to thank Anusch Taraz for telling her about the problem and some helpful discussions. 
Moreover, she expresses her thanks to Peter Heinig for useful hints that helped to make algorithms more efficient, which was very valuable in view of the numerical challenges. Thanks also to Michael Ritter, Christian Semrau, Alexandra Merz und Christoph Nie{\ss} for making the computers function and contributing important information about software.

\bibliographystyle{alpha}   
\bibliography{K_4paper.bib}  

\end{document}